\documentclass[review,3p,times,12pt]{elsarticle}

\usepackage{amssymb}
\usepackage{latexsym}
\usepackage{amsmath}
\usepackage{amssymb}
\usepackage{pifont}
\usepackage[colorlinks]{hyperref}
\usepackage{amsthm}
\newtheorem{theorem}{Theorem}[section]

\newtheorem{remark}[theorem]{Remark}
\newtheorem{example} [theorem]{Example}

\newdefinition{defn}[thm]{Definition}
\newdefinition{rmk}[thm]{Remark}
\newdefinition{ex}[thm]{Example}
\newproof{pf}{Proof}

\journal{Approximation Theory and Analytic Inequalities volume, edited by Th.M. Rassias, Springer}

\begin{document}

\begin{frontmatter}



\title{Hyers-Ulam  stability for differential equations and partial differential equations via Gronwall Lemma}


\author[sc]{Sorina Anamaria Ciplea}
\ead{sorina.ciplea@ccm.utcluj.ro}
\author[dm]{Daniela Marian\corref{cor1}}
\ead{daniela.marian@math.utcluj.ro}
\author[dm]{Nicolaie Lungu}
\ead{nlungu@math.utcluj.ro}
\author[tr]{Themistocles M. Rassias}
\ead{trassias@math.ntua.gr}

\cortext[cor1]{Corresponding author.}
\address[dm]{Technical University of Cluj-Napoca, Department of Mathematics, 28 Memorandumului Street, 400114,
Cluj-Napoca, Romania}
\address[sc]{Technical University of Cluj-Napoca, Department of Management and Technology, 28 Memorandumului Street, 400114,
Cluj-Napoca, Romania}
\address[tr]{National Technical University of Athens, Department of Mathematics, Zografou Campus, 15780,
Athens,  Greece}

\begin{abstract}
In this paper we will study Hyers-Ulam  stability for Bernoulli differential equations, Riccati differential equations and quasilinear partial differential equations of first order, using Gronwall Lemma, following a method given by Rus.

\end{abstract}


\end{frontmatter}

\section{Introduction}
In \cite{Rus1}, \cite{Rus2}, \cite{Rus3}  Rus has obtained some results regarding Ulam stability  of differential and integral equations, using Gronwall inequalities method and  weak Picard operators technique. In \cite{RL} Rus and Lungu have studied the stability of a  partial differential equation of order two of hyperbolic type using the same method. In \cite{LC}  Craciun and Lungu have studied, using this method, a partial differential equation of order two having a general form. In this paper we use the same method in order to study the stability of Bernoulli and Riccati equations and also of quasilinear partial differential equations of first order.
We mention that some results regarding Ulam stability of Bernoulli and Riccati differential equations was established by Jung and Rassias \cite {JR1}, \cite{JR2}, using the integrating factor method.
The first result proved on the Hyers-Ulam stability of partial differential equations
is due to A. Prastaro and Th.M. Rassias \cite{PraRassias}.
Also Lungu and Popa \cite{LP} and Marian and Lungu \cite{LM} have obtained stability results from some partial differential linear and quasilinear equations. The Gronwall inequality is used in Quarawani \cite{Qar} in order to study Hyers-Ulam-Rassias stability for Bernoulli differential equations and it is also used in \cite{CL}, \cite{LR}. For a broader study of Hyers-Ulam stability for functional equations the reader is also referred to the following books and papers: \cite{Abdollahpour}, \cite{Brzdek}, \cite{Jung}, \cite{Jung2}, \cite{Jung3}, \cite{Jung4}, \cite{JR1}, \cite{JR2}, \cite{Kannappan},  \cite{Lee}, \cite{Lee2}, \cite{Milovanovic}, \cite{Mortici}, \cite{Park}, \cite{Ulam}.

In the following we will use Definition 2.1, 2.2, 2.3 from \cite{Rus1}, p.126 and Remark 2.1, 2.2. from \cite{Rus1}, p.127.

\section{Main results}
\subsection{Stability of Bernoulli differential equation}

Let $\left(\mathbb{B},\left\vert \cdot \right\vert \right) $ be a (real or complex)
Banach space,  $a,b \in \mathbb{R}, a<b$, $p,q\in C\left( \left[ a,b\right] ,\mathbb{B}\right) $ and $n\in
\mathbb{R} \backslash \{0,1\}.$

We consider the Bernoulli  differential equation
\begin{equation}\label{eq:10}
 z^{\prime }\left( x\right) =p\left( x\right) z\left( x\right)
+q\left( x\right) z^{n}\left( x\right) ,
x\in \left[ a,b\right],
\end{equation}
and the inequation
\begin{equation}\label{ineq:10}
\left\vert y^{\prime }\left( x\right) -p\left( x\right) y\left( x\right)
-q\left( x\right) y^{n}\left( x\right) \right\vert \leq \varepsilon ,
x\in \left[ a,b\right] .
\end{equation}
From  Remark 2.1 from \cite{Rus1}, p.127 follows that $y\in C^{1}\left( \left[ a,b\right] ,\mathbb{B}\right) $
is a solution of the inequation \eqref{ineq:10} if and only if there exists a function $%
g\in C^{1}\left( \left[ a,b\right] ,\mathbb{B}\right) $ (which depend on $y$) such
that
\begin{itemize}
\item[(i)]  $\left\vert g\left( x\right) \right\vert \leq \varepsilon ,\forall x\in %
\left[ a,b\right] ;$

\item[(ii)] $y^{\prime }\left( x\right) =p\left( x\right) y\left( x\right) +q\left(
x\right) y^{n}\left( x\right) +g\left( x\right) ,\forall x\in \left[
a,b\right].$
\end{itemize}
From  Remark 2.2 from \cite{Rus1}, p.127 follows that if $y\in C^{1}\left( \left[ a,b\right] ,\mathbb{B}\right) $
is a solution of the inequation \eqref{ineq:10}, then $y$ is a solution of the following
integral inequation%
\[
\left\vert y\left( x\right) -y\left( a\right) -\int_{a}^{x}\left[
p\left( t\right) y\left( t\right) +q\left( t\right) y^{n}\left( t\right) %
\right] dt\right\vert \leq \left( x-a\right) \varepsilon ,\forall x\in \left[
a,b\right] .
\]

\begin{theorem}\label{th1}
If
\begin{itemize}
\item[(i)]  $a<\infty ,b<\infty ;$

\item[(ii)] $p,q\in C\left( \left[ a,b\right] ,\mathbb{B}\right) ;$

\item[(iii)] there exists $L>0$ such that%
\[
\left\vert q \left(x \right) y^{n}\left( x\right) - q \left(x \right) z^{n}\left( x\right) \right\vert \leq
L \left\vert y\left( x\right) -z\left( x\right) \right\vert ,
\]%
for all $x\in \left[ a,b\right] $ and $y,z\in C^{1}\left( \left[ a,b\right]
,\mathbb{B}\right) ,$
\end{itemize}
 then the equation \eqref{eq:10} is Hyers-Ulam  stable.
\end{theorem}

{\bf Proof.}{
Let $y\in C^{1}\left( \left[ a,b\right] ,\mathbb{B}\right) $ be a solution of the
inequation \eqref{ineq:10} and  $z$  the unique solution of the Cauchy problem
\begin{equation}
\left\{
\begin{array}{ll}
z^{\prime }\left( x\right)  =p\left( x\right) z\left( x\right) +q\left(
x\right) z^{n}\left( x\right) , x\in \left[ a,b\right] , \\
z\left( a\right) =y\left( a\right)
\end{array}
.\right.
\end{equation}
We have that%
\[
z\left( x\right) =y\left( a\right) +\int_{a}^{x}\left[ p\left(
t\right) z\left( t\right) +q\left( t\right) z^{n}\left( t\right) \right]
dt, x\in \left[ a,b\right] .
\]
Let
\begin{equation*}
M=\underset{x\in \left[ a,b\right] }{\max } \left \vert p\left( x\right) \right \vert .
\end{equation*}

We consider the difference
\begin{align*}
&\left\vert y\left( x\right) -z\left( x\right) \right\vert  \leq
\left\vert y\left( x\right) -y\left( a\right) -\int_{a}^{x}\left[ p\left(
t\right) y\left( t\right) +q\left( t\right) y^{n}\left( t\right) \right]
dt\right\vert + \\
&\left\vert \int_{a}^{x}\left[ p\left( t\right) y\left( t\right)
+q\left( t\right) y^{n}\left( t\right) -p\left( t\right) z\left( t\right)
-q\left( t\right) z^{n}\left( t\right) \right] dt\right\vert  \leq  \\
&\leq \varepsilon \left( x-a\right) +\int_{a}^{x}\left[ \left\vert p\left(
t\right)  y\left( t\right) - p\left(
t\right) z\left( t\right)
\right\vert +\left\vert \left(
 q\left( t\right) y^n\left( t\right) - q\left( t\right) z^n\left( t\right) \right) \right\vert \right] dt\leq\\
&\leq \varepsilon \left( x-a\right) +\int_{a}^{x}\left[ \left\vert p\left( t\right)\right\vert \left\vert y\left(t \right) -z\left( t\right) \right\vert +L \left\vert y\left(t \right) -z\left( t\right) \right\vert\right]dt=\\
&=\varepsilon \left( x-a\right) +\int_{a}^{x}\left[ \left\vert p\left( t\right)\right\vert+L\right] \left\vert y\left(t \right) -z\left( t\right) \right\vert dt
\end{align*}
From Gronwall lemma  (see \cite{LLM}, p. 6)  we have that
\begin{align*}
\left\vert y\left( x\right) -z\left( x\right) \right\vert  &\leq
  \varepsilon \left( x-a\right)
e^{\int_{a}^{x}\left[ \left\vert p\left(
t\right) \right\vert +L \right]dt}\leq \varepsilon \left( b-a\right)
e^{\int_{a}^{b}\left( M
+L\right )dt}=\\& = \varepsilon \left( b-a\right) e^{\left( M+L \right)\left( b-a \right)}= c \cdot \varepsilon,
\end{align*}
where $ c =   \left( b-a\right) e^{\left( M+L \right)\left( b-a \right)}.$
}

\begin{example}
We consider the Bernoulli differential equation%
\begin{equation}\label{eq:15}
z^{\prime }=xz+\frac{x}{1+x^{2}}\sqrt{z},
\end{equation}%
where $x\in \left[ a,b\right] $ and $z\geq 1.$We have $p\left( x\right) =x$
and $q\left( x\right) =\frac{x}{1+x^{2}}.$ Let $%
D=\left\{ \left( x,z\right) \mid x\in \left[ a,b\right] ,z\geq 1\right\} $
and $f\left( x,z\right) =\frac{x}{1+x^{2}}\sqrt{z}.$ We have%
\begin{equation*}
\left\vert \frac{\partial f}{\partial z}\right\vert =\left\vert \frac{x}{%
1+x^{2}}\cdot \frac{1}{2\sqrt{z}}\right\vert \leq \frac{1}{2}\left\vert
\frac{x}{1+x^{2}}\right\vert \leq \frac{1}{4},\forall \left( x,z\right) \in D,
\end{equation*}%
hence the function $f$ satisfies a Lipschitz condition in the variable $z$,
on $D$, with Lipschitz constant $1/4$. Hence%
\begin{equation*}
\left\vert f\left( x,y\right) -f\left( x,z\right) \right\vert \leq
L\left\vert y-z\right\vert =\frac{1}{4}\left\vert y-z\right\vert ,
\end{equation*}%
that is%
\begin{equation*}
\left\vert \frac{x}{1+x^{2}}\sqrt{y}-\frac{x}{1+x^{2}}\sqrt{z}\right\vert
\leq \frac{1}{4}\left\vert y-z\right\vert, x \in \left[a, b \right], y \geq 1, z \geq 1.
\end{equation*}
We apply Theorem \ref{th1} so the equation \eqref{eq:15} is Hyers-Ulam  stable.
Let $y\in C^{1}\left( \left[ a,b\right], \mathbb{B}\right) $ be a solution of the
inequation
\begin{equation}\label{eq:16}
\left\vert z^{\prime }-xz-\frac{x}{1+x^{2}}\sqrt{z}\right\vert \leq \varepsilon ,
\end{equation}%
and  $z$  the unique solution of the Cauchy problem
\begin{equation}
\left\{
\begin{array}{ll}
z^{\prime }=xz+\frac{x}{1+x^{2}}\sqrt{z},\\
z\left( a\right) =y\left( a\right)
\end{array}
.\right.
\end{equation}
We have
\[
z\left( x\right) =y\left( a\right) -\int_{a}^{x}\left[ t z+\frac{t}{1+t^{2}}\sqrt{z} \right]
dt, x\in \left[ a,b\right] .
\]
Let
\begin{equation*}
M=\underset{x\in \left[ a,b\right] }{\max } \left \vert p\left( x\right) \right \vert= \left \vert b \right \vert.
\end{equation*}
We have
\begin{equation*}
\left\vert y\left( x\right) -z\left( x\right) \right\vert \leq \varepsilon
\left( b-a\right) e^{\left( \left\vert b\right\vert +\frac{1}{4}  \right)
\left( b-a\right) }.
\end{equation*}
\end{example}

\bigskip

\subsection{Stability of Riccati differential equation}

Let $\left(\mathbb{B},\left\vert \cdot \right\vert \right) $ be a (real or complex)
Banach space,  $a,b \in \mathbb{R}, a<b$ and  $p,q,r\in C\left( \left[ a,b\right] ,\mathbb{B}\right).$

We consider the Riccati differential equation
\begin{equation}\label{eq:20}
 z^{\prime }\left( x\right) =p\left( x\right) z^2\left( x\right)
+q\left( x\right) z\left( x\right)+r\left( x\right),
x\in \left[ a,b\right],
\end{equation}
and the inequation
\begin{equation}\label{ineq:20}
\left\vert y^{\prime }\left( x\right) -p\left( x\right) y^2\left( x\right)
-q\left( x\right) y\left( x\right)-r\left( x \right) \right\vert \leq \varepsilon ,
x\in \left[ a,b\right].
\end{equation}
From Remark 2.1 from \cite{Rus1}, p.127  follows that $y\in C^{1}\left( \left[ a,b\right],\mathbb{B}\right) $
is a solution of the inequation \eqref{ineq:20} if and only if there exists a function $%
g\in C^{1}\left( \left[ a,b\right] ,\mathbb{B}\right) $ (which depend on $y$) such
that
\begin{itemize}
\item[(i)]  $\left\vert g\left( x\right) \right\vert \leq \varepsilon ,\forall x\in %
\left[ a,b\right];$

\item[(ii)] $y^{\prime }\left( x\right) =p\left( x\right) y^2\left( x\right) +q\left(
x\right) y\left( x\right)+r\left( x\right) +g\left( x\right) ,\forall x\in \left[
a,b\right].$
\end{itemize}

From  Remark 2.2 from \cite{Rus1}, p.127  follows that if $y\in C^{1}\left( \left[ a,b\right],\mathbb{B}\right) $
is a solution of the inequation \eqref{ineq:20}, then $y$ is a solution of the following
integral inequation%
\[
\left\vert y\left( x\right) -y\left( a\right) -\int_{a}^{x}\left[
p\left( t\right) y^2\left( t\right) +q\left( t\right) y\left( t\right) +r\left( t\right) %
\right] dt\right\vert \leq \left( x-a\right) \varepsilon ,\forall x\in \left[
a,b\right].
\]

\begin{theorem}
If
\begin{itemize}
\item[(i)]  $a<\infty ,b<\infty ;$

\item[(ii)] $p,q,r\in C\left( \left[ a,b\right] ,\mathbb{B}\right) ;$

\item[(iii)] there exists $L>0$ such that%
\[
\left\vert p \left(t \right) y^{2}\left( x\right) - p \left(t \right) z^{2}\left( x\right) \right\vert \leq
L \left\vert y\left( x\right) -z\left( x\right) \right\vert ,
\]%
for all $x\in \left[ a,b\right] $ and $y,z\in C^{1}\left( \left[ a,b\right]
,\mathbb{B}\right) ,$
\end{itemize}
 then the equation \eqref{eq:20} is Hyers-Ulam  stable.
\end{theorem}

{\bf Proof.}{
Let $y\in C^{1}\left( \left[ a,b\right],\mathbb{B}\right) $ be a solution of the
inequation \eqref{ineq:20} and  $z$  the unique solution of the Cauchy problem
\begin{equation}
\left\{
\begin{array}{ll}
z^{\prime }\left( x\right)  =p\left( x\right) z^2\left( x\right) +q\left(
x\right) z\left( x\right)+r\left( x\right) , x\in \left[ a,b\right], \\
z\left( a\right) =y\left( a\right)
\end{array}
.\right.
\end{equation}
We have that%
\[
z\left( x\right) =y\left( a\right) +\int_{a}^{x}\left[ p\left(
t\right) y^2\left( t\right) +q\left( t\right) z\left( t\right) +r\left( x\right) \right]
dt,\forall x\in \left[ a,b\right].
\]
Let
\begin{equation*}
M=\underset{x\in \left[ a,b\right] }{\max } \left \vert q\left( x\right) \right \vert .
\end{equation*}

We consider the difference
\begin{align*}
&\left\vert y\left( x\right) -z\left( x\right) \right\vert  \leq
 \left\vert y\left( x\right) -y\left( a\right) -\int_{a}^{x}\left[ p\left(
t\right) y^2\left( t\right) +q\left( t\right) y\left( t\right)+r\left( x\right) \right]
dt\right\vert + \\
&\left\vert \int_{a}^{x}\left[ p\left( t\right) y^2\left( t\right)
+q\left( t\right) y\left( t\right) -p\left( t\right) z^2\left( t\right)
-q\left( t\right) z\left( t\right) \right] dt\right\vert  \leq  \\
& \leq \varepsilon \left( x-a\right) +\int_{a}^{x} \left[ \left\vert p\left(
t\right)    y^2\left( t\right) -  p\left(
t\right) z^2\left( t\right)
\right\vert +\left\vert \left(
 q\left( t\right) y\left( t\right) - q\left( t\right) z\left( t\right) \right) \right\vert \right] dt \leq  \\
&\leq \varepsilon \left( x-a\right) +\int_{a}^{x} \left[ L \left\vert y\left( t\right) -z\left( t\right)
\right\vert +\left\vert q\left(
t\right) \right\vert  \left\vert y\left( t\right) -z\left( t\right) \right\vert
\right] dt = \\
&\leq \varepsilon \left( x-a\right) +\int_{a}^{x} \left[ L  +\left\vert q\left(
t\right) \right\vert  \right] \left\vert y\left( t\right) -z\left( t\right)
\right\vert dt .
\end{align*}
From Gronwall lemma  (see \cite{LLM}, p. 6)  we have that
\begin{align*}
&\left\vert y\left( x\right) -z\left( x\right) \right\vert  \leq
\varepsilon \left( x-a\right)
e^{\int_{a}^{x}\left(L +\left\vert q\left(
t\right) \right\vert  \right)dt}\leq \varepsilon \left( b-a\right)
e^{\int_{a}^{b}\left(L+M
 \right )dt}=\\
&= \varepsilon \left( b-a\right) e^{\left( M+L  \right)\left( b-a \right)}= c \cdot \varepsilon,
\end{align*}

where $ c= \left( b-a\right) e^{\left( M+L \right)\left( b-a \right)}.$
}

\bigskip

\subsection{Hyers-Ulam  stability of quasilinear partial differential equation}
\subsubsection{Hyers-Ulam  stability}

Let $\left(\mathbb{B},\left\vert \cdot \right\vert \right) $ be a (real or complex)
Banach space, $a,b\in \left( 0,\infty \right] $, $\varepsilon $  a positive
real number, $\varphi \in C\left( \left[ 0,a\right) \times \left[ 0,b\right)
,\mathbb{R}_{+}\right) $  and $p,q,r\in C\left( \left[ 0,a\right) \times \left[
0,b\right) \times \mathbb{B}, \mathbb{R}\right) $ and $p\left( x,y,u\right) \neq 0$.

We consider the following quasilinear partial differential equation of first
order%
\begin{equation}\label{eq:30}
\frac{\partial u\left( x,y\right) }{\partial x}=-\frac{q\left( x,y,u\right)
}{p\left( x,y,u\right) }\frac{\partial u}{\partial y}+\frac{r\left(
x,y,u\right) }{p\left( x,y,u\right) },
\end{equation}%
and the following partial differential inequation
\begin{equation}\label{eq:31}
\left\vert \frac{\partial v\left( x,y\right) }{\partial x}+\frac{q\left(
x,y,v\right) }{p\left( x,y,v\right) }\frac{\partial v}{\partial y}-\frac{%
r\left( x,y,v\right) }{p\left( x,y,v\right) }\right\vert \leq \varepsilon ,
\end{equation}

\begin{equation}\label{eq:32}
\left\vert \frac{\partial v\left( x,y\right) }{\partial x}+\frac{q\left(
x,y,v\right) }{p\left( x,y,v\right) }\frac{\partial v}{\partial y}-\frac{%
r\left( x,y,v\right) }{p\left( x,y,v\right) }\right\vert \leq \varepsilon
\cdot \varphi \left( x,y\right) .
\end{equation}

\begin{remark}\label{obs1}
A function $v\in C\left( \left[ 0,a\right) \times \left[ 0,b\right)
, \mathbb{B}\right) $ is a solution of the inequation \eqref{eq:31} if and only if there exists
a function $g\in C\left( \left[ 0,a\right) \times \left[ 0,b\right)
, \mathbb{B}\right) $  such that
\begin{itemize}
\item[(i)] $\left\vert g\left( x,y \right) \right\vert \leq \varepsilon ,\forall \left(
x,y\right) \in  \left[ 0,a\right) \times \left[ 0,b\right)$;

\item[(ii)] $\frac{\partial v\left( x,y\right) }{\partial x}=-\frac{q\left(
x,y,v\right) }{p\left( x,y,v\right) }v_{y}\left( x,y\right) +\frac{r\left(
x,y,v\right) }{p\left( x,y,v\right) }+g\left( x,y\right) $, where $v_{y}=%
\frac{\partial v}{\partial y}.$
\end{itemize}
\end{remark}

\begin{remark}\label{obs2}
If $v\in C\left( \left[ 0,a\right) \times \left[ 0,b\right) ,B\right) $ is a
solution of the inequation \eqref{eq:31}, then $v$ is a solution of the following
integral inequation
\begin{equation*}
\left\vert v\left( x,y\right) -v\left( 0,y\right) -\int_{0}^{x}\left[
-\frac{q\left( s,y,v\left( s,y\right) \right) }{p\left( s,y,v\left(
s,y\right) \right) }v_{y}\left( s,y\right) +\frac{r\left( s,y,v\left(
s,y\right) \right) }{p\left( x,y,v\left( s,y\right) \right) }\right]
ds\right\vert
 \leq \varepsilon x,
\end{equation*}
$\forall x\in \left[ 0,a\right) ,y\in \left[
0,b\right).$
\end{remark}

 Indeed, by Remark \ref{obs1} we have that
\begin{equation*}
\frac{\partial v\left( x,y\right) }{\partial x}=-\frac{q\left( x,y,v \left(x,y\right)\right)
}{p\left( x,y,v\left(x,y\right)\right) }v_{y}\left( x,y\right) +\frac{r\left( x,y,v \left(x,y\right)\right) }{%
p\left( x,y,v \left(x,y\right)\right) }+g\left( x,y\right),
\end{equation*}
$\forall x\in \left[ 0,a\right)
,y\in \left[ 0,b\right).$
This implies that
\begin{equation*}
v\left( x,y\right) =v\left( 0,y\right) +\int_{a}^{x}\left[ -\frac{%
q\left( s,y,v\left( s,y\right) \right) }{p\left( s,y,v\left( s,y\right)
\right) }v_{y}\left( s,y\right) +\frac{r\left( s,y,v\left( s,y\right)
\right) }{p\left( x,y,v\left( s,y\right) \right) }+g\left( s,y\right) \right]
ds.
\end{equation*}%
From this it follows that
\begin{align*}
&\left\vert v\left( x,y\right) -v\left( 0,y\right) -\int_{0}^{x}
\left[
-\frac{q\left( s,y,v\left( s,y\right) \right) }{p\left( s,y,v\left(
s,y\right) \right) }v_{y}\left( s,y\right) +\frac{r\left( s,y,v\left(
s,y\right) \right) }{p\left( x,y,v\left( s,y\right) \right) }\right]
ds\right\vert   \\
&\leq \int_{0}^{x}\left\vert g\left( s,y\right)
\right\vert ds\leq \varepsilon x.
\end{align*}

\begin{theorem}
We suppose that
\begin{itemize}
\item[(i)]  $a<\infty ,b<\infty ;$

\item[(ii)] $p,q,r\in C\left( \left[ 0,a\right] \times \left[ 0,b\right] \times
\mathbb{B}, \mathbb{B}\right) ,p\neq 0;$

\item[(iii)] there exists $l_{1},l_{2}>0$ such that%
\begin{align*}
&\left\vert \frac{q\left( x,y,v_{1}\right) }{p\left( x,y,v_{1}\right) }%
v_{1y}\left( x,y\right) -\frac{q\left( x,y,v_{2}\right) }{p\left(
x,y,v_{2}\right) }v_{2y}\left( x,y\right) \right\vert \leq l_{1}\left\vert
v_{1}-v_{2}\right\vert, \\
&\left\vert \frac{r\left( x,y,v_{1}\right) }{p\left( x,y,v_{1}\right) }-\frac{%
r\left( x,y,v_{2}\right) }{p\left( x,y,v_{2}\right) }\right\vert \leq
l_{2}\left\vert v_{1}-v_{2}\right\vert ,
\end{align*}%
$\forall v_{1},v_{2}\in \mathbb{B},\forall \left( x,y\right) \in \left[ 0,a\right]
\times \left[ 0,b\right].$
\end{itemize}
Then:
\begin{itemize}
\item[(a)] for $\psi \in C\left( \left[ 0,a\right] , \mathbb{B}\right) $ the equation \eqref{eq:30} has a
unique solution with
\begin{equation*}
u\left( 0,y\right) =\psi \left( y\right) ,\forall y\in \left[ 0,b\right] ;
\end{equation*}

\item[(b)] the equation \eqref{eq:30} is Hyers-Ulam  stable.
\end{itemize}
\end{theorem}

{\bf Proof.}{

(a) This is a known result (see \cite{Rus0} ).

(b) Let $v$ be a solution of the inequation \eqref{eq:31}. Denote by $u$ the unique
solution of the equation \eqref{eq:30} which satisfies the condition%
\begin{equation*}
u\left( 0,y\right) =v\left( 0,y\right)  ,\forall y\in \left[ 0,b\right] .
\end{equation*}

From Remark \ref{obs2} and condition (iii) we have that
\begin{align*}
&\left\vert v\left( x,y\right) -u\left( x,y\right) \right\vert  \leq\\
&\leq \left\vert v\left( x,y\right) -v\left( 0,y\right) -\int_{0}^{x}
\left[ -\frac{q\left( s,y,v\left( s,y\right) \right) }{p\left( s,y,v\left(
s,y\right) \right) }v_{y}\left( s,y\right) +\frac{r\left( s,y,v\left(
s,y\right) \right) }{p\left( x,y,v\left( s,y\right) \right) }\right]
ds\right\vert +\\
&\int_{0}^{x}\left\vert -\frac{q\left( s,y,v\left( s,y\right)
\right) }{p\left( s,y,v\left( s,y\right) \right) }v_{y}\left( s,y\right) +%
\frac{r\left( s,y,v\left( s,y\right) \right) }{p\left( x,y,v\left(
s,y\right) \right) }+\frac{q\left( s,y,u\left( s,y\right) \right) }{p\left(
s,y,u\left( s,y\right) \right) }u_{y}\left( s,y\right) -\frac{r\left(
s,y,u\left( s,y\right) \right) }{p\left( x,y,u\left( s,y\right) \right) }
\right\vert ds \\
&\leq \varepsilon x+\int_{0}^{x}\left[ l_{1}\left\vert v\left(
s,y\right) -u\left( s,y\right) \right\vert +l_{2}\left\vert v\left(
s,y\right) -u\left( s,y\right) \right\vert \right] ds.
\end{align*}
Or,
\begin{equation*}
\left\vert v\left( x,y\right) -u\left( x,y\right) \right\vert \leq
\varepsilon x+\int_{0}^{x}\left[ l_{1}+l_{2}\right] \left\vert
v\left( s,y\right) -u\left( s,y\right) \right\vert ds.
\end{equation*}

From Gronwall lemma (see \cite{LLM}, p. 6)  we have
\begin{equation*}
\left\vert v\left( x,y\right) -u\left( x,y\right) \right\vert \leq
ae^{a\left( l_{1}+l_{2}\right) }\cdot \varepsilon =c\cdot \varepsilon ,
\end{equation*}%
where $c=ae^{a\left( l_{1}+l_{2}\right) }.$

So, the equation \eqref{eq:30} is Hyers-Ulam  stable.
}
\subsubsection{Hyers-Ulam-Rassias stability of equation \eqref{eq:30}}

Let us consider the equation \eqref{eq:30} and the inequation \eqref{eq:32} in the case $a=\infty
,b=\infty .$

\begin{theorem}
We suppose that
\begin{itemize}
\item[(i)]  $p,q,r\in C\left( \left[ 0,a\right) \times \left[ 0,b\right) \times
\mathbb{B}, \mathbb{B}\right) ,p\neq 0;$

\item[(ii)] there exists $l_{1}, l_{2} \in C^{1}\left( \left[ 0, a\right) \times \left[ 0,b\right), \mathbb{R}_{+}\right)$ such that
\begin{align*}
&\left\vert \frac{q\left( x,y,v_{1}\right) }{p\left( x,y,v_{1}\right) }%
v_{1y}\left( x,y\right) -\frac{q\left( x,y,v_{2}\right) }{p\left(
x,y,v_{2}\right) }v_{2y}\left( x,y\right) \right\vert \leq l_{1}\left(
x,y\right) \left\vert v_{1}-v_{2}\right\vert;\\
&\left\vert \frac{r\left( x,y,v_{1}\right) }{p\left( x,y,v_{1}\right) }-\frac{%
r\left( x,y,v_{2}\right) }{p\left( x,y,v_{2}\right) }\right\vert \leq
l_{2}\left( x,y\right) \left\vert v_{1}-v_{2}\right\vert;
\end{align*}%
$\forall v_{1},v_{2}\in \mathbb{B},\forall \left( x,y\right) \in \left[ 0,a\right)
\times \left[ 0,b\right) ;$

\item[(iii)] $e^{\int_{0}^{\infty}\left[ l_{1}\left( s,y\right) +l_{2}\left( s,y\right) \right] ds}$ is convergent and there exists a real number $M$ such that
 $e^{\int_{0}^{\infty}\left[ l_{1}\left( s,y\right) +l_{2}\left( s,y\right) \right] ds}\leq M,$ $\forall y \in \left[ 0,b\right)$;

\item[(iv)] there exists $\lambda _{\varphi }>0$ such that%
\begin{equation*}
\int_{0}^{x}\varphi \left( s,y\right) ds\leq \lambda _{\varphi }
 \cdot \varphi \left( x,y\right) , \forall \left( x,y\right) \in \left[ 0,a\right)
\times \left[ 0,b\right)
\end{equation*}%
and $\varphi $ increasing.
\end{itemize}
Then the equation \eqref{eq:30} ($a=\infty ,b=\infty $) is
Hyers-Ulam-Rassias stable.
\end{theorem}

{\bf Proof.}{
Let $v$ be a solution of the inequation \eqref{eq:32}. Denote by $u$ the unique solution
of the problem%
\begin{equation*}
\left\{
\begin{array}
[c]{l}%
\frac{\partial u\left( x,y\right) }{\partial x}=-\frac{q\left( x,y,u\right)
}{p\left( x,y,u\right) }u_{y} \left( x,y\right)+\frac{r\left(
x,y,u\right) }{p\left( x,y,u\right) } \\
u\left( 0,y\right) =v\left( 0,y\right) .%
\end{array}%
\right.
\end{equation*}
We have%
\begin{equation*}
u\left( x,y\right) =v\left( 0,y\right) +\int_{0}^{x}\left[ -\frac{q\left(
s,y,u\left( s,y\right) \right) }{p\left( s,y,u\left( s,y\right) \right) }%
u_{y}\left( s,y\right) +\frac{r\left( s,y,u\left( s,y\right) \right) }{%
p\left( x,y,u\left( s,y\right) \right) }\right] ds
\end{equation*}%
and%
\begin{align*}
&\left\vert v\left( x,y\right) -v\left( 0,y\right) -\int_{0}^{x}\left[ -\frac{%
q\left( s,y,v\left( s,y\right) \right) }{p\left( s,y,v\left( s,y\right)
\right) }v_{y}\left( s,y\right) +\frac{r\left( s,y,v\left( s,y\right)
\right) }{p\left( x,y,v\left( s,y\right) \right) }\right] ds\right\vert \leq\\
&\leq \varepsilon \int_{0}^{x}\varphi \left( s,y\right) ds\leq \varepsilon \lambda
_{\varphi } \cdot \varphi \left( x,y\right).
\end{align*}%
Then we have
\begin{align*}
&\left\vert v\left( x,y\right) -u\left( x,y\right) \right\vert  \leq\\
&\leq \left\vert v\left( x,y\right) -v\left( 0,y\right) -\int_{0}^{x}\left[ -%
\frac{q\left( s,y,v\left( s,y\right) \right) }{p\left( s,y,v\left(
s,y\right) \right) }v_{y}\left( s,y\right) +\frac{r\left( s,y,v\left(
s,y\right) \right) }{p\left( x,y,v\left( s,y\right) \right) }\right]
ds\right\vert + \\
&\int_{0}^{x}\left\vert -\frac{q\left( s,y,v\left( s,y\right) \right)
}{p\left( s,y,v\left( s,y\right) \right) }v_{y}\left( s,y\right) +\frac{%
r\left( s,y,v\left( s,y\right) \right) }{p\left( x,y,v\left( s,y\right)
\right) }+\frac{q\left( s,y,u\left( s,y\right) \right) }{p\left( s,y,u\left(
s,y\right) \right) }u_{y}\left( s,y\right) -\frac{r\left( s,y,u\left(
s,y\right) \right) }{p\left( x,y,u\left( s,y\right) \right) }\right\vert
ds \\
&\leq \varepsilon \lambda _{\varphi } \cdot \varphi \left(
x,y\right) +\int_{0}^{x}\left[ l_{1}\left( s,y\right) \left\vert v\left( s,y\right)
-u\left( s,y\right) \right\vert +l_{2}\left( s,y\right) \left\vert v\left(
s,y\right) -u\left( s,y\right) \right\vert \right]ds \leq\\
&\leq \varepsilon \lambda _{\varphi } \cdot \varphi \left(
x,y\right) +\int_{0}^{x}\left[ l_{1}\left( s,y\right) +l_{2}\left(
s,y\right) \right] \left\vert v\left( s,y\right) -u\left( s,y\right)
\right\vert ds.
\end{align*}
From Gronwall lemma  (see \cite{LLM}, p. 6)   we have that
\begin{equation*}
\left\vert v\left( x,y\right) -u\left( x,y\right) \right\vert \leq
\varepsilon \lambda _{\varphi } \cdot \varphi \left(
x,y\right) e^{\int_{0}^{x}\left[ l_{1}\left( s,y\right) +l_{2}\left(
s,y\right) \right] ds}\leq c_{\varphi }\cdot \varepsilon \cdot \varphi \left( x,y\right) ,
\end{equation*}%
where $c_{\varphi }= \lambda _{\varphi } \cdot M.$

So, the equation \eqref{eq:30} is generalized Hyers-Ulam-Rassias stable.
}

\bibliographystyle{elsarticle-harv}

\end{document}